\documentclass[12pt]{article}
\usepackage{amsthm,amsmath,amssymb,amscd,verbatim}
\usepackage{graphics,pifont}
\usepackage{amssymb}
\usepackage{graphicx}
\def\sq{$\square$}

\date{}

\newcommand{\ora}{\overrightarrow}

\begin{document}
\title{On the Class of Similar Square $\{ -1, 0, 1 \}$-Matrices Arising from Vertex maps on Trees}
\author{Bau-Sen Du \cr
Institute of Mathematics \cr
Academia Sinica \cr
Taipei 10617, Taiwan \cr
dubs@math.sinica.edu.tw \cr}
\maketitle

\begin{abstract}
Let $n \ge 2$ be an integer.  In this note, we show that the {\it oriented} transition matrices over the field $\mathcal R$ of all real numbers (over the finite field $\mathcal Z_2$ of two elements respectively) of all continuous {\it vertex maps} on {\it all} oriented trees with $n+1$ vertices are similar to one another over $\mathcal R$ (over $\mathcal Z_2$ respectively) and have characteristic polynomial $\sum_{k=0}^n x^k$.  Consequently, the {\it unoriented} transition matrices over the field $Z_2$ of all continuous {\it vertex maps} on {\it all} oriented trees with $n+1$ vertices are similar to one another over $\mathcal Z_2$ and have characteristic polynomial $\sum_{k=0}^n x^k$.  Therefore, the coefficients of the characteristic polynomials of these {\it unoriented} transition matrices, when considered over the field $\mathcal R$, are all odd integers (and hence nonzero).    


\end{abstract}

Let $n \ge 2$ be an integer and let $T$ be a tree with $n+1$ vertices $V_1, V_2, \cdots, V_{n+1}$.  The tree $T$ has $n$ edges, say, $E_1, E_2, \cdots, E_n$.  If vertices $V_{j_1}$ and $V_{j_2}$ are endpoints of an edge $E$, then, we let $[V_{j_1} : V_{j_2}]$ denote the edge $E$, i.e., the set of all points in $E$ and, following Bernhardt {\bf\cite{b1}}, we denote the {\it positively oriented edge} from $V_{j_1}$ to $V_{j_2}$ as $\overrightarrow{[V_{j_1}, V_{j_2}]}$ and call $V_{j_1}$ the first vertex of $\overrightarrow{[V_{j_1}, V_{j_2}]}$ and $V_{j_2}$ the second.  We also define $-\overrightarrow{[V_{j_1}, V_{j_2}]}$ by putting $-\overrightarrow{[V_{j_1}, V_{j_2}]} = \overrightarrow{[V_{j_2}, V_{j_1}]}$ and call it the {\it negatively oriented edge} from $V_{j_1}$ to $V_{j_2}$.  So, the first vertex of $-\overrightarrow{[V_{j_1}, V_{j_2}]}$ is $V_{j_2}$ and the second is $V_{j_1}$.  Thus, both $\overrightarrow{[V_{j_1}, V_{j_2}]}$ and $-\overrightarrow{[V_{j_1}, V_{j_2}]}$ represent the same edge $[V_{j_1} : V_{j_2}]$, but with the opposite orientations.  In the sequel, we denote these $n$ positively oriented edges of $T$ as $\ora{E_1}, \ora{E_2}, \cdots, \ora{E_n}$ and call the resulting tree {\it oriented tree} and denote it as $\ora{T}$.  It is clear that there are exactly $2^n$ distinct such oriented trees $\overrightarrow{T}$.  Later, we shall see that different choices of orientations on the edges of $T$ will not affect our main results.  When no confusion arises, we shall always use $V_1, V_2, \cdots, V_{n+1}$ and $\ora{E_1}, \ora{E_2}, \cdots, \ora{E_n}$ to denote respectively the vertices and the oriented edges of {\it any} tree with $n+1$ vertices.

Following {\bf\cite{b1}}, for any two vertices $V_i$ and $V_j$ in the oriented tree $\ora{T}$, a path from $V_i$ to $V_j$ is a sequence of {\it oriented edges} $\ora{e_1}, \ora{e_2}, \cdots, \ora{e_m}$, where $\ora{e_k} \in \{ \ora{E_s}, -\ora{E_s} : 1 \le s \le n \}$ for all $1 \le k \le m$, and the first vertex of $\ora{e_1}$ is $V_i$, the second vertex of $\ora{e_m}$ is $V_j$, and the second vertex of $\ora{e_\ell}$ is equal to the first vertex of $\ora{e_{\ell+1}}$ for all $1 \le \ell \le m-1$.  It is clear that, for any two vertices, $V_i$ and $V_j$, there is a {\it unique shortest path} from $V_i$ to $V_j$ in $\ora{T}$ which will be denoted as $\ora{[V_i, V_j]}$.  We also denote $-\ora{[V_i, V_j]} = \ora{[V_j, V_i]}$ as the shortest path from vertex $V_j$ to vertex $V_i$.  From now on, when we write {\it the shortest path} ${\ora{[V_i, V_j]}}$, we always mean the {\it shortest path} from vertex $V_i$ to vertex $V_j$ in $\ora{T}$.  We also let $[V_i : V_j]$ denote the collection of all points in all (oriented) edges in {\it the shortest path} $\ora{[V_i, V_j]}$.  

Following {\bf\cite{b1, b2}}, let $f : T \longrightarrow T$ be a continuous {\it vertex map}, i.e., $f$ is a continuous map such that the $n+1$ vertices of $T$ form a periodic orbit and, for each $1 \le i \le n$, $f$ is {\it monotonic} on the {\it unoriented edge} $E_i = [V_{i_1} : V_{i_2}]$, meaning that, as the point $x$ moves from vertex $V_{i_1}$ to vertex $V_{i_2}$ {\it monotonically} along the edge $E_i$, the point $f(x)$ moves {\it monotonically} from vertex $f(V_{i_1})$ to vertex $f(V_{i_2})$ along {\it the shortest path} $\overrightarrow{[f(V_{i_1}), f(V_{i_2})]}$ from $f(V_{i_1})$ to $f(V_{i_2})$.  We now let $\mathcal R$ denote the field of all real numbers and let $\mathcal Z_2$ denote the finite field $\{ 0, 1 \}$ of two elements and let $\mathcal F$ be any field with unity 1.  We define the associated {\it oriented transition} $n \times n$ $\{ -1, 0, 1 \}$-matrix $\mathcal A_n(f) = (\alpha_{i,j})$ over $\mathcal F$ by putting the {\it positively oriented edge} $\ora{E_i} = \ora{[V_{i_1}, V_{i_2}]}$ and putting 
$$
\alpha_{i,j} = \begin{cases}
               1, & \text{if $\ora{E_j}$ appears in {\it the shortest path} $\ora{[f(V_{i_1}), f(V_{i_2})]}$ \small{from vertex} $f(V_{i_1})$ \small{to vertex} $f(V_{i_2})$}, \cr
               -1, & \text{if $-\ora{E_j}$ appears in {\it the shortest path} $\ora{[f(V_{i_1}), f(V_{i_2})]}$ \small{from vertex} $f(V_{i_1})$ \small{to vertex} $f(V_{i_2})$}, \cr
               0, & \text{otherwise}. \cr
       \end{cases}
$$
and define the associated {\it unoriented transition} $n \times n$ $\{ 0, 1 \}$-matrix $\mathcal B_n(f) = (\beta_{i,j})$ over $\mathcal F$ by putting, for all $1 \le i \le n$ and all $1 \le j \le n$, $\beta_{i,j} = 1$ if $\alpha_{i,j} \ne 0$ and $\beta_{i,j} = 0$ otherwise, or equivalently, 
$$
\beta_{i,j} = \begin{cases}
               1, & \text{if the set inclusion $f(E_i) \supset E_j$ holds}, \cr
               0, & \text{otherwise}. \cr
       \end{cases}
$$
There are exactly $2^n$ such oriented transition matrices $\mathcal A_n(f)$ for each $f$ and yet they all have the same unoriented transition matrix $\mathcal B_n(f)$.  Later, we shall see that the determinant of $\mathcal A_n(f)$ is $(-1)^n$ and that of $\mathcal B_n(f)$ is an odd integer.  In {\bf\cite{du1}}, we study the special case when $T$ is a compact interval in the real line and $f$ is a continuous vertex map on $T$.  In this note, we generalize the main results in {\bf\cite{du1}} for interval maps to vertex maps on trees.  Surprisingly, the arguments used there {\it almost} work for vertex maps on trees.  For completeness, we include the proofs.  

Recall that $\mathcal F$ denotes a field with unity 1.  Let $\ora {W_{\mathcal F}}^n(\mathcal E) = \bigl\{ \sum_{i=1}^{n} r_i\ora{E_i} : r_i \in \mathcal F, 1 \le i \le n \bigr\}$ denote the $n$-dimensional vector space over $\mathcal F$ with $\mathcal E = \{ \ora{E_j} : 1 \le j \le n \}$ as a basis.  In the sequel, when there is no confusion, we shall write $\ora {W_{\mathcal F}}^n$ instead of $\ora {W_{\mathcal F}}^n(\mathcal E)$.  So, now we regard each positively oriented edge $\ora{E_j}$ as a basis element of the vector space $\ora {W_{\mathcal F}}^n$ and regard the negatively oriented edge $-\ora{E_j}$ of $\ora{E_j}$ as an element in $\ora {W_{\mathcal F}}^n$ such that $\ora{E_j} + (-\ora{E_j}) = {\bf 0}$.  Let $\sum_{i=1}^n r_i\ora{E_i}$ be an element of $\ora{W_{\mathcal F}}^n$ such that $r_i \in \{ -1, 0, 1 \}$ for all $1 \le i \le n$.  If there exist two vertices $V_i$ and $V_j$ such that $r_k = 1$ if and only if $\ora{E_k}$ appears in {\it the shortest path} $\ora{[V_i, V_j]}$ from veretx $V_i$ to vertex $V_j$ and $r_k = -1$ if and only if $-\ora{E_k}$ appears in {\it the shortest path} $\ora{[V_i, V_j]}$ from vertex $V_i$ to vertex $V_j$, then we define $\ora{[V_i, V_j]} = \sum_{i=1}^n r_i\ora{E_i}$ and $\ora{[V_j, V_i]} = -\ora{[V_i, V_j]} = -(\sum_{i=1}^n r_i\ora{E_i})$.  In particular, if $\ora{[V_i, V_j]} = \ora{e_1}\ora{e_2} \cdots \ora{e_m}$ is {\it the shortest path} in the oriented tree $\ora{T}$ defined as above, then, as elements of $\ora{W_{\mathcal F}}^n$, we have $\ora{[V_i, V_j]} = \sum_{k=1}^m \ora{e_k}$.  Therefore, for any two vertices $V_i$ and $V_j$, the notation $\ora{[V_i, V_j]}$ will have two meanings: It represents the {\it unique shortest path} from vertex $V_i$ to vertex $V_j$ in the oriented tree $\ora{T}$ on the one hand, and represents the element of $\ora{W_{\mathcal F}}^n$ which is a sum of those oriented (positively or negatively) edges which appear in the {\it unique shortest path} $\overrightarrow{[V_i, V_j]}$ from vertex $V_i$ to vertex $V_j$ on the other.  There should be no confusion from the texts.  With respect to the {\it oriented transition} $n \times n$ $\{ -1, 0, 1 \}$-matrices $\mathcal A_n(f) = (\alpha_{i,j})$ of the continuous vertex tree map $f$, we define a linear transformation $\Phi_f$ from $\ora{W_{\mathcal F}}^n$ into itself such that, for each $1 \le i \le n$, $$\Phi_f(\ora{E_i}) = \sum_{j=1}^n \alpha_{i,j} \ora{E_j}.$$Therefore, if $\ora{E_i} = \ora{[V_{i_1}, V_{i_2}]}$ is a positively oriented edge of $\ora{T}$ from vertex $V_{i_1}$ to vertex $V_{i_2}$, then, when considered as an element of $\ora{W_{\mathcal F}}^n$, we have, by definition of $\mathcal A_n(f)$ and $\Phi_f$, $\Phi_f(\ora{E_i}) = \sum_{j=1}^n \alpha_{ij}\ora{E_j} = \ora{[f(V_{i_1}), f(V_{i_2})]}$ {\small which also represents the} {\it unique shortest path} {\small from vertex} $f(V_{i_1})$ {\small to vertex} $f(V_{i_2})$ in $\ora{T}$.  

We shall need the following fundamental result.

\noindent
{\bf Lemma 1.}
{\it For any distinct vertices $V_i$ and $V_j$ of the oriented tree $\ora{T}$, we have $\Phi_f(\ora{[V_i, V_j]}) = \ora{[f(V_i), f(V_j)]}$.  That is, if $\ora{[V_i, V_j]}$ is the unique shortest path from vertex $V_i$ to vertex $V_j$ in $\ora{T}$, then $\Phi_f(\ora{[V_i, V_j]}) \, (= \ora{[f(V_i), f(V_j)]})$ is the unique shortest path from vertex $f(V_i)$ to vertex $f(V_j)$ in $\ora{T}$.  Similarly, if $V_{i_1}, V_{i_2}, \cdots, V_{i_m}$ are vertices of $\ora{T}$, then $\sum_{k=1}^{m-1} \ora{[V_{i_k}, V_{i_{k+1}}]} = \ora{[V_{i_1}, V_{i_m}]}$.}  

\noindent
{\it Proof.}
Let $V_i, V_k, V_j$ be three distinct vertices of the tree $\ora{T}$.  Assume that both $\ora{[V_i, V_k]}$ and $\ora{[V_k, V_j]}$ are positively oriented edges of $\ora{T}$.  If the set intersection $f([V_i : V_k]) \cap f([V_k : V_j]) = \{ f(V_k) \}$, then the concatenation of {\it the shortest path} $\ora{[f(V_i), f(V_k)]}$ and {\it the shortest path} $\ora{[f(V_k), f(V_j)]}$ becomes the {\it shortest} path $\ora{[f(V_i), f(V_j)]}$ from vertex $f(V_i)$ to vertex $f(V_j)$.  Therefore, we have $\Phi_f(\ora{[V_i, V_j]}) = \ora{[f(V_i), f(V_j)]}$.  On the other hand, if the set intersection $f([V_i : V_k]) \cap f([V_k : V_j]) = [V_\ell : f(V_k)] \ne  \{ f(V_k) \}$ for some veterx $V_\ell \ne f(V_k)$, then {\it the shortest path} $\ora{[V_\ell, f(V_k)]}$ in {\it the shortest path} $\ora{[f(V_i), f(V_k)]}$ and {\it the shortest path} $\ora{[f(V_k), V_\ell]}$ in {\it the shortest path} $\ora{[f(V_k), f(V_j)]}$ cancel out.  So, $\Phi_f(\ora{[V_i, V_j]}) = \Phi_f(\ora{[V_i, V_k]} + \ora{[V_k, V_j]}) = \Phi_f(\ora{[V_i, V_k]}) + \Phi_f(\ora{[V_k, V_j]}) = \ora{[f(V_i), f(V_k)]} + \ora{[f(V_k), f(V_j)]} =\bigr(\ora{[f(V_i), V_\ell]} + \ora{[V_\ell, f(V_k)]}\bigr) + \bigr(\ora{[f(V_k), V_\ell]} + \ora{[V_\ell, f(V_j)]}\bigr) = \ora{[f(V_i), V_\ell]} + \ora{[V_\ell, f(V_j)]} = \ora{[f(V_i), f(V_j)]}$.   

Assume that both $\ora{[V_i, V_k]}$ and $\ora{[V_j, V_k]} \,(= -\ora{[V_k, V_j]})$ are positively oriented edges of the oriented tree $\ora{T}$.  Then {\it the shortest path} $\ora{[V_i, V_j]}$ is the concatenation of the positively oriented edge $\ora{[V_i, V_k]}$ and the negatively oriented edge $(-\ora{[V_j, V_k]})$.  Thus, as elements of $\ora{W_{\mathcal F}}^n$, we have $\ora{[V_i, V_j]} = \ora{[V_i, V_k]} - \ora{[V_j, V_k]}$.  So, $\Phi_f(\ora{[V_i, V_j]}) = \Phi_f(\ora{[V_i, V_k]} - \ora{[V_j, V_k]}) = \Phi_f(\ora{[V_i, V_k]}) - \Phi_f(\ora{[V_j, V_k]}) = \ora{[f(V_i), f(V_k)]} - \ora{[f(V_j), f(V_k)]}$.  If the set intersection $f([V_i : V_k]) \cap f([V_j : V_k]) = \{ f(V_k) \}$, then the concatenation of {\it the shortest path} $\ora{[f(V_i), f(V_k)]}$ and {\it the shortest path} $-\ora{[f(V_j), f(V_k)]} \, (= \ora{[f(V_k), f(V_j)]})$ becomes the {\it shortest} path $\ora{[f(V_i), f(V_j)]}$ from vertex $f(V_i)$ to vertex $f(V_j)$.  Therefore, we obtain that $\Phi_f(\ora{[V_i, V_j]}) = \ora{[f(V_i), f(V_j)]}$.  On the other hand, if the set intersection $f([V_i : V_k]) \cap f([V_j : V_k]) = [V_\ell : f(V_k)] \ne  \{ f(V_k) \}$ for some veterx $V_\ell \ne f(V_k)$, then {\it the shortest path} $\ora{[V_\ell, f(V_k)]}$ in {\it the shortest path} $\ora{[f(V_i), f(V_k)]}$ and {\it the shortest path} $-\ora{[V_\ell, f(V_k)]} \, (=\ora{[f(V_k), V_\ell]})$ in {\it the shortest path} $-\ora{[f(V_j), f(V_k)]}) \, (=\ora{[f(V_k), f(V_j)]}$) cancel out.  Therefore, $\Phi_f(\ora{[V_i, V_j]}) = \Phi_f(\ora{[V_i, V_k]} - \ora{[V_j, V_k]}) = \Phi_f(\ora{[V_i, V_k]}) - \Phi_f(\ora{[V_j, V_k]}) = \ora{[f(V_i), f(V_k)]} - \ora{[f(V_j), f(V_k)]} =\bigr(\ora{[f(V_i), V_\ell]} + \ora{[V_\ell, f(V_k)]}\bigr) - \bigr(\ora{[f(V_j), V_\ell]} + \ora{[V_\ell, f(V_k)]}\bigr) = \ora{[f(V_i), V_\ell]} - \ora{[f(V_j), V_\ell]} = \ora{[f(V_i), V_\ell]} + \ora{[V_\ell, f(V_j)]} = \ora{[f(V_i), f(V_j)]}$.

If both $\overrightarrow{[V_k, V_i]} \, (= -\overrightarrow{[V_i, V_k]})$ and $\overrightarrow{[V_k, V_j]}$ or, both $\overrightarrow{[V_k, V_i]} \, (= -\overrightarrow{[V_i, V_k]})$ and $\overrightarrow{[V_j, V_k]} \, (= -\overrightarrow{[V_k, V_j]}$ are positively oriented edges of the tree $\ora{T}$, then, by discussing cases depending on the set intersections $f([V_i : V_k]) \cap f([V_j : V_k])$ as above, we obtain that $\Phi_f(\ora{[V_i, V_j]}) = \ora{[f(V_i), f(V_j)]}$.  We omit the details.  

So far, we have shown that $\Phi_f(\ora{[V_i, V_j]}) = \ora{[f(V_i), f(V_j)]}$ as long as {\it the shortest path} $\ora{[V_i, V_j]}$ consists of exactly two oriented edges.  Now, if {\it the shortest path} $\ora{[V_i, V_j]} = \ora{e_1}\ora{e_2}\ora{e_3}$ consists of exactly three oriented edges $\ora{e_1}, \ora{e_2}, \ora{e_3}$.  Let the second vertex of $\ora{e_2}$ be $V_k$.  Then $\ora{[V_i, V_k]} = \ora{e_1}\ora{e_2}$.  It follows from what we just proved above that $\Phi_f(\ora{[V_i, V_k]}) = \ora{[f(V_i), f(V_k)]}$.  Therefore, $\Phi_f(\ora{[V_i, V_j]}) = \Phi_f(\ora{e_1}\ora{e_2}\ora{e_3}) = \Phi_f(\ora{e_1}\ora{e_2}+\ora{e_3}) = \Phi_f(\ora{e_1}\ora{e_2}) + \Phi_f(\ora{e_3}) = \Phi_f(\ora{[V_i, V_k]}) + \Phi_f(\ora{e_3})$.  Depending on $\ora{e_3} = \ora{[V_k, V_j]}$ or $\ora{e_3} = -\ora{[V_k, V_j]}$ and arguing as above, we can easily obtain that $\Phi_f(\ora{[V_i, V_j]}) = \ora{[f(V_i), f(V_j)]}$ whenever $\ora{[V_i, V_j]} = \ora{e_1}\ora{e_2}\ora{e_3}$ consists of exactly three oriented edges $\ora{e_1}, \ora{e_2}, \ora{e_3}$.  The general case when $\ora{[V_i, V_j]}$ consists of more than 3 oriented edges can be proved similarly by induction.  Therefore, $\Phi_f(\ora{[V_i, V_j]}) = \ora{[f(V_i), f(V_j)]}$ as long as $V_i$ and $V_j$ are any two distinct vertices of $\ora{T}$.  

Finally, if $V_{i_1}, V_{i_2}, \cdots, V_{i_m}$ are $m \ge 2$ distinct vertices of $\ora{T}$, then similar arguments show that $\sum_{k=1}^{m-1} \ora{[V_{i_k}, V_{i_{k+1}}]} = \ora{[V_{i_1}, V_{i_m}]}$.  This completes the proof.
\hfill\sq

\noindent
{\bf Lemma 2.}
{\it $\Phi_f$ is an isomorphism from $\ora{W_{\mathcal F}}^n$ onto itself.}

\noindent
{\it Proof.}
Let $\hat f$ be any continuous vertex map on the tree $T$ such that the composition $\hat f \circ f$ is the identity map on the vertices of $T$.  Then, by Lemma 1, for each positively oriented edge $\ora{E_i} = \ora{[V_{i_1}, V_{i_2}]}$, we have $(\Phi_{\hat f} \circ \Phi_f)({\ora{E_i}}) = \Phi_{\hat f}(\Phi_f(\ora{[V_{i_1}, V_{i_2}]})) = \Phi_{\hat f}(\ora{[f(V_{i_1}), f(V_{i_2})]}) = \ora{[(\hat f \circ f)(V_{i_1}), (\hat f \circ f)(V_{i_2})]} = \ora{[V_{i_1}, V_{i_2}]}$ $= \ora{E_i}$.  Therefore, $\Phi_{\hat f}$ is the inverse of $\Phi_f$.
\hfill\sq

We shall need the following result which is proved in {\bf\cite{du1}}.  For completeness, we include its proof.  

\noindent
{\bf Lemma 3.}
{\it Let $1 \le j \le n$ be any fixed integer and let $b$ denote the greatest common divisor of $j$ and $n+1$.  Let $s = (n+1)/b$.  For every integer $1 \le k \le s-1$, let $1 \le m_k \le n$ be the unique integer such that $kj \equiv m_k$ (mod $n+1$).  Then the $m_k$'s are all distinct and $\{ m_k : 1 \le k \le s-1 \} = \{ kb : 1 \le k \le s-1 \}$.}

\noindent
{\it Proof.}
Let $B = \{ m_k : 1 \le k \le s-1 \}$ and $C = \{ kb : 1 \le k \le s-1 \}$.  For every integer $1 \le k \le s-1$, since $j/b$ and $(n+1)/b$ are relatively prime, the congruence equation $(j/b)x \equiv k$ (mod $(n+1)/b$) has an integer solution $x$ such that $1 \le x \le s-1 = [(n+1)/b] -1$.  Consequently, for every integer $1 \le k \le s-1$, the congruence equation $(m_x \equiv) \, jx \equiv kb$ (mod $n+1$) has an integer solution $x$ such that $1 \le x \le s-1$.  Since $1 \le kb \le (s-1)b \le n$ and $1 \le m_k \le n$ for every integer $1 \le k \le s-1$, we obtain that $C \subset B$.  Since both $B$ and $C$ consist of exactly $s-1$ elements, we obtain that $B = C$.  That is, $\{ m_k : 1 \le k \le s-1 \} = \{ kb : 1 \le k \le s-1 \}$.  This completes the proof.
\hfill\sq

Let $M_1$ and $M_2$ be two $n \times n$ matrices over the field $\mathcal F$.  We say that $M_1$ is similar to $M_2$ through the invertible matrix $G$ if $M_1 \cdot G = G \cdot M_2$.  We can now prove our main result.

\noindent
{\bf Theorem 1.}
{\it Let $n \ge 2$ be an integer.  Let $T$ be any tree with $n+1$ vertices.  Let $f$ be a continuous vertex map on $T$.  Let $\mathcal R$, $\mathcal Z_2$, $\mathcal F$, $\ora{W_{\mathcal F}}^n$, $\ora{W_{\mathcal Z_2}}^n$, $\Phi_f$, $\mathcal A_n(f)$ and $\mathcal B_n(f)$ be defined as above.  Then the following hold:
\begin{itemize}
\item[{\rm (1)}]
For each integer $1 \le i \le n$, $\sum_{k=1}^n \Phi_f^k(\ora{E_i}) = {\bf 0}$ and so, $\sum_{k=1}^n \Phi_f^k(w) = {\bf 0}$ for all $w$ in $\ora{W_{\mathcal F}}^n$.  

\item[{\rm (2)}]
Let $i$ and $j$ be two integers in the interval $[1, n]$ and let $\ora{J}$ denote {\it the shortest path} $\ora{[V_i, f^j(V_i)]}$ in $\ora{T}$.  If $j$ and $n+1$ are relatively prime, then the set $\mathcal W_f = \{ \Phi_f^k(\ora{J}) : 0 \le k \le n-1 \}$ is a basis for $\ora{W_{\mathcal Z_2}}^n$ and for $\ora{W_{\mathcal F}}^n$ when $\mathcal F$ is a field with characteristic zero or the determinant of the matrix $\mathcal M_f$ of the set $\mathcal W_f$ with respect to the basis $\mathcal E = \{ \ora{E_1}, \ora{E_2}, \cdots, \ora{E_n} \}$ is not "divisible" by the finite characteristic of $\mathcal F$.  Furthermore, when $T$ is a tree in the real line with $n+1$ vertices and $f$ is a continuous vertex map on $T$, then the constant term of the characteristic polynomial of the matrix $\mathcal M_f$ is $\pm 1$ and hence the set $\mathcal W_f$ is a basis of $\ora{W_{\mathcal F}}^n$ for any field $\mathcal F$ (however, not all coefficients of the characteristic polynomial of the matrix $\mathcal M_f$ are odd integers (see, for example, Figure 1(a) with $\ora{J} = \ora{[1, 2]}$ where the corresponding characteristic polynomial is $x^5-x^4-x+1$).

\item[{\rm (3)}]
Over any field $\mathcal F$ with characteristic zero ($\mathcal Z_2$ respectively), the oriented transition matrix $\mathcal A_n(f)$ and its inverse $[\mathcal A_n(f)]^{-1}$, as $\{ -1, 0, 1 \}$-matrices, are similar to the following companion matrix 
$$
\left[ {\begin{array}{*{20}c} 0 & 1 &  0 &  0 &  \cdots & 0  \\ 0 &  0 &  1 &  0 & \cdots & 0  \\ 0 & 0 &  0 &  1 & \cdots & 0 \\ & \vdots &   &   & \vdots &  \\  -1 &  -1 &  -1 & -1 &  \cdots & -1  \\  \end{array}} \right]
$$
of the polynomial $\sum_{k=0}^{n} x^k$ through invertible $\{ -1, 0, 1 \}$-matrices over $\mathcal F$ ($\mathcal Z_2$ respectively) and have the same characteristic polynomial $\sum_{k=0}^{n} x^k$ while the unoriented transition matrices of all continuous vertex maps on $T$ may not be similar to each other over the same field $\mathcal F$ ($\mathcal Z_2$ respectively) (see Figures 1 - 3).  Furthermore, if $T$ is a tree in the real line, then the oriented transition matrices, when considered over any field, of all continuous vertex maps on $T$ with $n+1$ vertices and their inverses are similar to one another through invertible $\{ -1, 0, 1 \}$-matrices and have the same characteristic polynomial $\sum_{k=0}^{n} x^k$.
 
\item[{\rm (4)}]
The coefficients of the characteristic polynomial of the unoriented transition matrix $\mathcal B_n(f)$, when considered as a matrix over $\mathcal R$, are all odd integers (see Figures 1 - 4).  Furthermore, the unoriented transition matrices over any field of all continuous vertex maps on all trees with $n+1$ vertices, when considered as matrices over $\mathcal Z_2$, are similar to one another and have characteristic polynomial $\sum_{k=0}^{n} x^k$, but may not be similar to each other when considered over the finite field $\mathcal Z_p = \{ 0, 1, 2, \cdots, p-1 \}$, where $p \ge 3$ is a prime number (see Figures 1 - 3).    
\end{itemize}}

\noindent
{\it Proof.}
To prove Part (1), recall that $f$ is a continuous vertex map on the tree $T_1$.  For any fixed integer $1 \le i \le n$, let $\ora{E_i} = \ora{[V_{i_1}, V_{i_2}]}$ and let $1 \le j \le n$ be the unique integer such that $f^j(\ora{V_{i_1}}) = \ora{V_{i_2}}$.  So, $\ora{E_i} = \ora{[V_{i_1}, V_{i_2}]} = \ora{[V_{i_1}, f^j(V_{i_1})]}$.  Let $b$ be the greatest common divisor of $j$ and $n+1$ and let $s = (n+1)/b$.  So, $sj = (j/b)(sb) = (j/b)(n+1)$.  For every integer $1 \le k \le s-1$, let $1 \le m_k \le n$ be the unique integer such that $kj \equiv m_k$ (mod $n+1$).  Then, by Lemma 3, we obtain that $\{ m_k : 1 \le k \le s-1 \} = \{ kb : 1 \le k \le s-1 \}$.  Let $m_0 = 0$.  Then $\{ m_k : 0 \le k \le s-1 \} = \{ kb : 0 \le k \le s-1 \}$.  Hence, the set $\{ 1, 2, \cdots, n-1, n \}$ is the disjoint union of the sets $\{ m_k + m : 0 \le k \le s-1 \}, 0 \le m \le b-1$.  Therefore, by Lemma 1, $\sum_{k=0}^{s-1} \Phi_f^{m_k}(\ora{E_i}) = \sum_{k=0}^{s-1} \Phi_f^{kj}(\ora{E_i})$ (since $kj \equiv m_k$ (mod $n+1)) = \ora{[V_{i_1}, f^j(V_{i_1})]} + \ora{[f^j(V_{i_1}), f^{2j}(V_{i_1})]} + \ora{[f^{2j}(V_{i_1}), f^{3j}(V_{i_1})]} + \cdots + \ora{[f^{(s-2)j}(V_{i_1}), f^{(s-1)j}(V_{i_1})]} + \ora{[f^{(s-1)j}(V_{i_1}), V_{i_1}]} = {\bf 0}$.  Thus, $\sum_{\ell=0}^n \Phi_f^\ell(\ora{E_i}) = \sum_{m=0}^{b-1} \Phi_f^m\bigr(\sum_{k=0}^{s-1} \Phi_f^{m_k}(\ora{E_i})\bigr) = {\bf 0}$.  Therefore, $\sum_{k=1}^n \Phi_f^k(w) = {\bf 0}$ for all vectors $w$ in $\ora{W_{\mathcal F}}^n$.  This establishes Part (1).  

For the proof of Part (2), we first consider those $\Phi_f$ over the field $\mathcal Z_2$.  Now we want to show that if $\mathcal N$ is a nonempty subset of $\{ 1, 2, \cdots, n-1, n \}$ such that $\ora{J} + \sum_{k \in \mathcal N} \Phi_f^k(\ora{J}) = \bf 0$, then $\mathcal N = \{ 1, 2, \cdots, n-1, n \}$.  Indeed, for every integer $1 \le k \le n$, let $1 \le m_k \le n$ be the unique integer such that $kj \equiv m_k$ (mod $n+1$).  Assume that $(j =) \,\, m_1 \notin \mathcal N$.  Then, for any $m \in \mathcal N$, $m \ne j$.  So, $f^m(V_i) \ne f^j(V_i)$.  If $(f^m(f^j(V_i)) =) \,\,f^{m+j}(V_i) = f^j(V_i)$, then the least period of $f^j(V_i)$ under $f$ divides $m$ ($< n+1$) which contradicts the fact that its least period under $f$ is $n+1$.  Therefore, {\it the shortest path} $\Phi_f^m(\ora{J}) = \Phi_f^m(\ora{[V_i, f^j(V_i)]}) = \ora{[f^m(V_i), f^{m+j}(V_i)]}$ either contains the vertex $f^j(V_i)$ in its "interior" or does not contain it.  So, in the expression of the element $\Phi_f^m(\ora{J}) = \Phi_f^m(\ora{[V_i, f^j(V_i)]}) = \ora{[f^m(V_i), f^{m+j}(V_i)]}$ as a sum of the basis elements $\ora{E_k}$'s, the number of the basis elements $\ora{E_k}$ which contain the vertex $f^j(V_i)$ as an endpoint is either 0 or 2.  Since $\ora{J} = \ora{[V_i, f^j(V_i)]}$ contains exactly one baisis element $\ora{E_k}$ which has the vertex $f^j(V_i)$ as an endpoint, there are an {\it odd} number of basis elements $\ora{E_k}$'s which has the vertex $f^j(V_i)$ as an endpoint in the expression of the element $\ora{J} + \sum_{k \in \mathcal N} \Phi_f^k(\ora{J})$ as a sum of the basis elements $\ora{E_k}$'s.  Consequently, $\ora{J} + \sum_{k \in \mathcal N} \Phi_f^k(\ora{J}) \ne {\bf 0}$.  This is a contradiction.  So, $(j =) \,\, m_1 \in \mathcal N$ and 
$$
{\bf 0} = \ora{J} + \sum_{k \in \mathcal N} \Phi_f^k(\ora{J}) = \ora{J} + \Phi_f^{m_1}(\ora{J}) + \sum_{k \in \mathcal N \setminus \{ m_1 \}} \Phi_f^k(\ora{J})
$$
$$
\,\, = \ora{[V_i, f^j(V_i)]} + \ora{[f^j(V_i), f^{2j}(V_i)]} + \sum_{k \in \mathcal N \setminus \{ m_1 \}} \Phi_f^k(J)$$ $$= \ora{[V_i, f^{2j}(V_i)]} + \sum_{k \in \mathcal N \setminus \{ m_1 \}} \Phi_f^k(\ora{J}).\qquad \qquad\qquad\quad
$$
Now assume that $m_2 \notin \mathcal N \setminus \{ m_1 \}$.  Then, for any $m \in \mathcal N \setminus \{ m_1 \}$, $m \notin \{ m_1, m_2 \} = \{ j, m_2 \} \subset \{ 1, 2, \cdots, n \}$.  If $f^m(V_i) = f^{2j}(V_i) \, (= f^{m_2}(V_i))$, then $m \equiv m_2$ (mod $n+1$).  Since both $m$ and $m_2$ are integers in the set $\{ 1, 2, \cdots, n \}$ such that $m \equiv m_2$ (mod $n+1$), we have $m = m_2$.  This is a contradiction.  If $f^{m+j}(V_i) = f^{2j}(V_i)$, then $m+j \equiv 2j$ (mod $n+1$) and so, $m \equiv j \, (\equiv m_1)$ (mod $n+1$).  Since both $m$ and $m_1$ are integers in the set $\{ 1, 2, \cdots, n \}$ such that $m \equiv m_1$ (mod $n+1$), we have $m = m_1$.  This is again a contradiction.  Therefore, in the expression of the element $\Phi_f^m(\ora{J}) = \Phi_f(\ora{[V_i, f^j(V_i)]}) = \ora{[f^m(V_i), f^{m+j}(V_i)]}$ as a sum of the basis elements $\ora{E_k}$'s, the number of the basis elements $E_k$ which contain the vertex $f^{2j}(V_i)$ as an endpoint is either 0 or 2.  Since $\ora{[V_i, f^{2j}(V_i)]}$ contains exactly one baisis element $\ora{E_k}$ which has the vertex $f^{2j}(V_i)$ as an endpoint, there are an {\it odd} number of basis elements $\ora{E_k}$'s which has the vertex $f^{2j}(V_i)$ as an endpoint in the expression of the element $\ora{[V_i, f^{2j}(V_i)]} + \sum_{k \in \mathcal N \setminus \{ m_1 \}} \Phi_f^k(\ora{J})$ as a sum of the basis elements $\ora{E_k}$'s.  Consequently, $\ora{J} + \sum_{k \in \mathcal N} \Phi_f^k(\ora{J}) = \ora{[V_i, f^{2j}(V_i)]} + \sum_{k \in \mathcal N \setminus \{ m_1 \}} \Phi_f^k(\ora{J})\ne {\bf 0}$.  This is a contradiction.  So, $m_2 \in \mathcal N \setminus \{ m_1 \}$ and $$
{\bf 0} = \ora{J} + \sum_{k \in \mathcal N} \Phi_f^k(\ora{J}) = \ora{[V_i, f^{2j}(V_i)]} + \sum_{k \in \mathcal N \setminus \{ m_1 \}} \Phi_f^k(\ora{J}).\qquad\qquad\qquad\quad\qquad\qquad\qquad\quad
$$
$$
=\ora{[V_i, f^{2j}(V_i)]} + \Phi_f^{m_2}(\ora{J}) + \sum_{k \in \mathcal N \setminus \{ m_1, m_2 \}} \Phi_f^k(\ora{J})\,\,
$$
$$
\qquad\quad\,\, = \ora{[V_i, f^{2j}(V_i)]} + \ora{[f^{2j}(V_i), f^{3j}(V_i)]} + \sum_{k \in \mathcal N \setminus \{ m_1, m_2 \}} \Phi_f^k(\ora{J})
$$
$$
\quad\quad\,\,= \ora{[V_i, f^{3j}(V_i)]} + \sum_{k \in \mathcal N \setminus \{ m_1, m_2 \}} \Phi_f^k(\ora{J}).\qquad\qquad\qquad\quad
$$
Proceeding in this manner finitely many times, we obtain that $\{ m_1, m_2, \cdots, m_{n-1}, m_n \} \subset \mathcal N$.  Since $j$ and $n+1$ are relatively prime, we see that, by Lemma 3, $\{ m_1, m_2, \cdots, m_n \} = \{ 1, 2, \cdots, n-1, n \}$.  Since $\{ m_1, m_2, \cdots, m_n \} \subset \mathcal N \subset \{ 1, 2, \cdots, n-1, n \}$, we obtain that $\mathcal N = \{ 1, 2, \cdots, n-1, n \}$.  This proves our assertion.  

Now assume that $\sum_{k=0}^{n-1} r_k\Phi_f^k(\ora{J}) = {\bf 0}$, where $r_k = 0$ or 1 in $\mathcal Z_2$, for all $0 \le k \le n-1$.  If $r_0 = 0$ and $r_\ell \ne 0$ for some integer $1 \le \ell \le n-1$, we may assume that $\ell$ is the smallest such integer.  Since $\Phi_f$ is invertible on $\ora{W_{\mathcal Z_2}}^n$, we obtain that $\ora{J} + \sum_{k=1}^{n-1-\ell} r_k\Phi_f^k(\ora{J}) = {\bf 0}$.  So, without loss of generality, we may assume that $r_0 \ne 0$.  That is, we may assume that $\ora{J} + \sum_{k=1}^{n-1} r_k\Phi_f^k(\ora{J}) = {\bf 0}$.  Let $\mathcal N = \{ k : 1 \le k \le n-1$ and $r_k \ne 0 \}$.  Then we have $\ora{J} + \sum_{k \in \mathcal N} \Phi_f^k(\ora{J}) = {\bf 0}$.  However, it follows from what we just proved above that $\mathcal N = \{ 1, 2, \cdots, n-1, n \}$.  This contradicts the assumption that $\mathcal N \subset \{ 1, 2, \cdots, n-1 \}$.  Therefore, the set $\{ \Phi_f^k(\ora{J}) : 0 \le k \le n-1 \}$ is linearly independent in the $n$-dimensional vector space $\ora{W_{\mathcal Z_2}}^n$ and hence is a basis for $\ora{W_{\mathcal Z_2}}^n$.  Consequently, the matrix of the basis $\{ \Phi_f^k(\ora{J}) : 0 \le k \le n-1 \}$ over $\mathcal Z_2$ with respect to the basis $\{ \ora{E_1}, \ora{E_2}, \cdots, \ora{E_n} \}$ of $\ora{W_{\mathcal Z_2}}^n$, denoted as 
$$
\bigl [\ora{J}, \Phi_f(\ora{J}), \Phi_f^2(\ora{J}), \cdots, \Phi_f^{n-1}(\ora{J}) : \ora{E_1}, \ora{E_2}, \cdots, \ora{E_n} \bigr ] \,\,\, (\text{over} \,\,\, \mathcal Z_2),
$$
has nonzero determinant and hence equals 1.  This implies that, over the general field $\mathcal F$ with unity 1, the $\{ -1, 0, 1 \}$-matrix 
$$
\mathcal M_f = \bigl [\ora{J}, \Phi_f(\ora{J}), \Phi_f^2(\ora{J}), \cdots, \Phi_f^{n-1}(\ora{J}) : \ora{E_1}, \ora{E_2}, \cdots, \ora{E_n} \bigr ] \,\,\, (\text{over} \,\,\, \mathcal F)
$$
of the set $\{ \Phi_f^k(\ora{J}) : 0 \le k \le n-1 \}$ with respect to the basis $\{ \ora{E_1}, \ora{E_2}, \cdots, \ora{E_n} \}$ of $\ora{W_{\mathcal F}}^n$ also has nonzero determinant if the characteristic of $\mathcal F$ is zero or the determinant of $\mathcal M_f$ is not {\it divisible} by the finite characteristic of $\mathcal F$.  In particular, if $\mathcal F = \mathcal R$ or $\mathcal F = \mathcal Z_2$, then the set $\{ \Phi_f^k(\ora{J}) : 0 \le k \le n-1 \}$ is a basis for $\ora{W_{\mathcal F}}^n$ {\bf\cite{her}}.  Note that when $\mathcal F = \mathcal R$, the determinant of $\mathcal M_f$ is an odd integer.  We do not know if it is always equal to $\pm 1$.  However, when $T$ is a tree in the real line, by choosing all orientations on the edges {\it same direction}, the matrix $\mathcal M_f$ over $\mathcal R$ is a Petrie matrix, i.e., in any row, nonzero entries are {\it consecutive} and are all equal to 1 or to $-1$.  It follows from easy induction {\bf\cite{gor}} that the determinant of any Petrie matrix is 0 or $\pm 1$.  Since the determinant of $\mathcal M_f$ over $\mathcal R$ is nonzero, we obtain that the determinant of $\mathcal M_f$ is $\pm 1$.  The rest is easy and omitted.  This confirms Part (2).

Now, for the general field $\mathcal F$ with unity 1, let $\ora{T}$ be the oriented tree on the interval $[1, n+1]$ in the real line with $n+1$ vertices $\hat V_i = i, 1 \le i \le n+1$ and $n$ positively oriented edges $\ora{D_j} = \ora{[j, j+1]}, 1 \le j \le n$.  Let $h$ be the continuous vertex map on $\ora{T}$ such that $h(x) = x+1$ for all $1 \le x \le n$ and $h(x) = -nx + n^2 + n + 1$ for all $n \le x \le n+1$.  Then $\Phi_h(\ora{D_k}) = \ora{D_{k+1}} = \Phi_h^{k}(\ora{D_1})$ for all $1 \le k \le n-1$ and $\Phi_h(\ora{D_n}) = \Phi_h(\ora{[\hat V_n, \hat V_{n+1}]}) = \ora{[\hat V_{n+1}, \hat V_1]} = -\ora{[\hat V_1, \hat V_{n+1}]} = -\sum_{k=1}^{n} \ora{D_k}$.  By definition, the set $\mathcal D = \{ \Phi_h^k(\ora{D_1}) : 0 \le k \le n-1 \} = \{ \ora{D_1}, \ora{D_2}, \cdots, \ora{D_n} \}$ is a basis for $\ora{W_{\mathcal F}}^n(\mathcal D)$.  Let $1 \le i \le n$ be a fixed integer and choose a fixed integer $1 \le j \le n$ such that $j$ and $n+1$ are relatively prime and let $\ora{J} = \ora{[V_i, f^j(V_i)]}$.  
$$
\text{Suppose the set} \,\,\, \{ \Phi_f^k(\ora{J}) : 0 \le k \le n-1 \} \,\,\, \text{is a basis for} \,\,\, \ora{W_{\mathcal F}}^n(\mathcal E). \qquad\qquad\qquad\qquad (*)
$$
Let $\phi : \ora{W_{\mathcal F}}^n(\mathcal D) \longrightarrow \ora{W_{\mathcal F}}^n(\mathcal E)$ be the linear transformation defined by $$\phi(\ora{D_k}) = \Phi_f^{k-1}(\ora{J}) \,\, \text{for all} \,\, 1 \le k \le n.$$Then $\phi$ is an isomorphism and the matrix of the basis $\{ \Phi_f^k(\ora{J}) : 0 \le k \le n-1 \}$ with respect to the basis $\mathcal E = \{ \ora{E_1}, \ora{E_2}, \cdots, \ora{E_n} \}$ of $\ora{W_{\mathcal F}}^n(\mathcal E)$ is an $n \times n$ $\{ -1, 0, 1 \}$-matrix.  Furthermore, 
$$
(\phi \circ \Phi_h)(\ora{D_n}) = \phi(\Phi_h(\ora{D_n})) = \phi(-\sum_{k=1}^{n} \ora{D_k}) = -\sum_{k=1}^{n} \phi(\ora{D_k}) = -\sum_{k=1}^{n} \Phi_f^{k-1}(\ora{J}) = \Phi_f^n(\ora{J}) \,\, \text{(by Part (1))}
$$
$$
\,\,\,\,= \Phi_f(\Phi_f^{n-1}(\ora{J})) = \Phi_f(\phi(\ora{D_n})) = (\Phi_f \circ \phi)(\ora{D_n})
$$
and, for every integer $1 \le k \le n-1$, 
$$
(\phi \circ \Phi_h)(\ora{D_k}) = \phi(\Phi_h(\ora{D_k})) = \phi(\ora{D_{k+1}}) = \Phi_f^k(\ora{J}) = \Phi_f(\Phi_f^{k-1}(\ora{J})) = \Phi_f(\phi(\ora{D_k})) = (\Phi_f \circ \phi)(\ora{D_k}).
$$
Therefore, $\Phi_f$ is similar to $\Phi_h$ through $\phi$.  Similarly, $\Phi_{\hat f}$ is similar to $\Phi_h$, where $\hat f$ is any continuous {\it vertex map} on the oriented tree $\ora{T}$ such that the composition $\hat f \circ f$ on the vertices of $\ora{T}$ is the identity map.  So, the matrices $\mathcal A_n(f)$ and $\mathcal A_n(\hat f) = [\mathcal A_n(f)]^{-1}$ are similar to $\mathcal A_n(h)$ over $\mathcal F$.  Similarly, the matrices $\mathcal A_n(g)$ and $[\mathcal A_n(g)]^{-1}$ are similar to $\mathcal A_n(h)$ over $\mathcal F$.  Consequently, we obtain that the matrices $\mathcal A_n(f)$, $[\mathcal A_n(f)]^{-1}$, $\mathcal A_n(g)$, and $[\mathcal A_n(g)]^{-1}$ are similar to one another over $\mathcal F$.  By Part (2), the above (*) holds for $\mathcal F = \mathcal Z_2$ and for any field $\mathcal F$ with char($\mathcal F) = 0$.  Therefore, the matrices $\mathcal A_n(f)$, $[\mathcal A_n(f)]^{-1}$, $\mathcal A_n(g)$, and $[\mathcal A_n(g)]^{-1}$ are similar to one another over $\mathcal Z_2$ and over any field  $\mathcal F$ with char($\mathcal F) = 0$.  On the other hand, let $P_n(x) = x^n + \cdots$ denote the characteristic polynomial of $\mathcal A_n(f)$ over $\mathcal Z_2$ or over a field $\mathcal F$ with char($\mathcal F) = 0$.  By Part (2), the degree of the {\it minimal} polynomial of the element $\ora{[V_i, f(V_i)]}$ is at least $n$.  It follows from Part (1) that the polynomial $\sum_{k=0}^{n} x^k$ is the minimal polynomial of $\ora{[V_i, f(V_i)]}$.  By the well-known Cayley-Hamilton theorem on matrices, we see that the element $\ora{[V_i, f(V_i)]}$ also satisfies the polynomial $P_n(x) - \sum_{k=0}^{n} x^k$ whose degree is at most $n-1 \, (< n)$.  Therefore, $P_n(x) - \sum_{k=0}^{n} x^k = 0$, i.e., the characteristic polynomial of $\mathcal A_n(f)$ is $\sum_{k=0}^{n} x^k$.  This proves Part (3).  

Let $\mathcal A_n(f)$ and $\mathcal B_n(f)$ be the oriented and unoriented transition matrices of $\Phi_f$ over $\mathcal R$ respectively.  Then it follows from Part (3) that the characteristic polynomial of $\mathcal A_n(f)$ is $\sum_{k=0}^n x^k$.  When we consider $\mathcal A_n(f)$ as a matrix over $\mathcal Z_2$, we obtain that $\mathcal B_n(f) = \mathcal A_n(f)$ and the characteristic polynomial of $\mathcal B_n(f)$ is $\sum_{k=0}^n x^k$ over $\mathcal Z_2$.  Consequently, the coefficients of the characteristic polynomial of $\mathcal B_n(f)$ over $\mathcal R$ are all odd integers (see Figures 1 - 4).  Furthermore, we see that $\mathcal B_n(f) = \mathcal A_n(f)$ over $\mathcal Z_2$.  So, it follows from Part (3) that $\mathcal B_n(f)$ and $\mathcal B_n(g)$, when considered as matrices over $\mathcal Z_2$, are similar to each other.  This proves Part (4) and completes the proof of the theorem.
\hfill\sq

\noindent
{\bf Remark.}
Let $f$ be a continuous vertex map on the tree $T \,\, (= T_1)$ with $n+1 \ge 3$ vertices.  For $\mathcal F = \mathcal R$ and any choices of orientations on the edges of $T$, it follows from Theorem 1(1) that the determinant of the corresponding oriented transition matrix $\mathcal A_n(f)$ is $(-1)^n$ while, by Theorem 1(4), that of the corresponding unoriented transition matrix $\mathcal B_n(f)$ is an odd integer which is not necessarily equal to $\pm 1$.  See Figure 2 for some examples.  In the following, we present two sufficient conditions which guarantee that the determinant of the corresponding unoriented transition matrix $\mathcal B_n(f)$ is $\pm 1$.  For other related problems regarding the unoriented transition matrices $\mathcal B_n(f)$, we refer to {\bf\cite{du2}} where (new) notions of one-sided and two-sided similarities or weak similarities of square $\{ 0, 1 \}$-matrices are introduced and examples are presented.  It is clear that notions of various similarities of the unoriented transition matrices similar to those considered in {\bf\cite{du2}} can be generalized from trees in the real line (i.e., compact intervals) to arbitrary trees.  

\noindent
{\bf Proposition 1.}
{\it Let $f$ be a continuous vertex map on the tree $T \,\, (= T_1)$ with $n+1 \ge 3$ vertices.  For each integer $1 \le i \le n$, let $\ora{E_i} = \ora{[V_{i_1}, V_{i_2}]}$ be a positively oriented edge of the oriented tree $\ora{T}$.  Let the field $\mathcal F = \mathcal R$.  Then, by Lemma 1 and the definition of the map $\Phi_f$ on the vector space $\ora{W_{\mathcal R}}^n$, we obtain that $\Phi_f(\ora{E_i}) = \ora{[f(V_{i_1}), f(V_{i_2})]}$.  Since $f$ is a continuous vertex map on the connected edge $E_i$ of the tree $T$, we can write $\Phi_f(\ora{E_i}) = \ora{[f(V_{i_1}), f(V_{i_2})]} = \sum_{j=1}^{m_i} r_{i,j} \ora{[V_{\ell_{i,j}}, V_{\ell_{i,j+1}}]}$, where, for each $1 \le j \le m_i$, $r_{i,j} = \pm 1$, $V_{\ell_{i,1}} = f(V_{i_1})$, $V_{\ell_{i,m_i+1}} = f(V_{i_2})$ and $\ora{[V_{\ell_{i,j}}, V_{\ell_{i,j+1}}]}$ is a positively oriented edge of the oriented tree $\ora{T}$.  Then the following hold:
\begin{itemize}
\item[{\rm (1)}]
If, for each $1 \le i \le n$, $\Phi_f(\ora{E_i})$ has only one sign, i.e., $r_{i,1} = r_{i,2} = \cdots = r_{i, m_i}$ (this includes the cases when $T$ is a compact interval in the real line), then the corresponding oriented transition matrix $\mathcal A_n(f)$ of $f$ can be obtained from that of the corresponding unoriented transition matrix $\mathcal B_n(f)$ of $f$ by performing the following row operation: Multiplying one row by $-1$.  Consequently, the determinant of the matrix $\mathcal B_n(f)$ is equal to $\pm 1$ times that of the matrix $\mathcal A_n(f)$ which is $\pm 1$ (see {\bf\cite{du1}}).  

\item[{\rm (2)}]
If, for each $1 \le i \le n$ such that $\Phi_f(\ora{E_i})$ does not have one sign, there exists an integer $1 \le k_i < m_i$ such that $r_{i,1} = r_{i,2} = \cdots = r_{i, k_i} \ne r_{i, k_i+1} = r_{i, k_i+2} = \cdots = r_{i, m_i}$ and $|r_{i, 1}| = |r_{i, m_i}| = 1$, let $\hat V_{\ell_{i,k_i}}$ be the unique vertex of $T$ such that $f(\hat V_{\ell_{i,k_i}}) = V_{\ell_{i,k_i}}$ and let $\ora{e_1} \ora {e_2} \cdots \ora{e_s}$ be the {\rm shortest} path from either $V_{i_1}$ or $V_{i_2}$ to $\hat V_{\ell_{i,k_i}}$ which passes through the edge $E_i = [V_{i_1}, V_{i_2}]$ (and so the second vertex of $\ora {e_s}$ is $\hat V_{\ell_{i,k_i}}$).  If each one of $\Phi_f(\ora{e_2}), \Phi_f(\ora{e_3}), \cdots, \Phi_f(\ora{e_{s-1}})$ and $\Phi_f(\ora{e_s})$ has only one sign, then the corresponding oriented transition matrix $\mathcal A_n(f)$ of $f$ can be obtained from that of the corresponding unoriented transition matrix $\mathcal B_n(f)$ of $f$ by performing the following two row operations: (i) Multiplying one row by $-1$ and (ii) Multiplying one row by $\pm 2$ and adding to another row.  Consequently, the determinant of the matrix $\mathcal B_n(f)$ is equal to $\pm 1$ times that of the matrix $\mathcal A_n(f)$ which is $\pm 1$ (see Figures 1 $\&$ 3).  
\end{itemize}}

\noindent
{\bf Remark.}
Figure 4 demonstrates a case which is not covered by Proposition 1, yet has the same conclusion.  We note that although, for a continuous vertex map $f$ on a tree $T$ with $n+1$ vertices, there are $2^n$ distinct oriented transition matrices, they all have one and the same unoriented transition matrix.  Therefore, if we can find an orientation for the tree $T$ so that Proposition 1 applies, then we obtain that the determinant of the unoriented transition matrix is $\pm 1$.  Figure 4 is such an example.

\bigskip
\begin{figure}[htbp]  
\centering
\begin{minipage}[t]{0.23\linewidth}
\centering
\includegraphics[width=1.6in, height=1.5in]{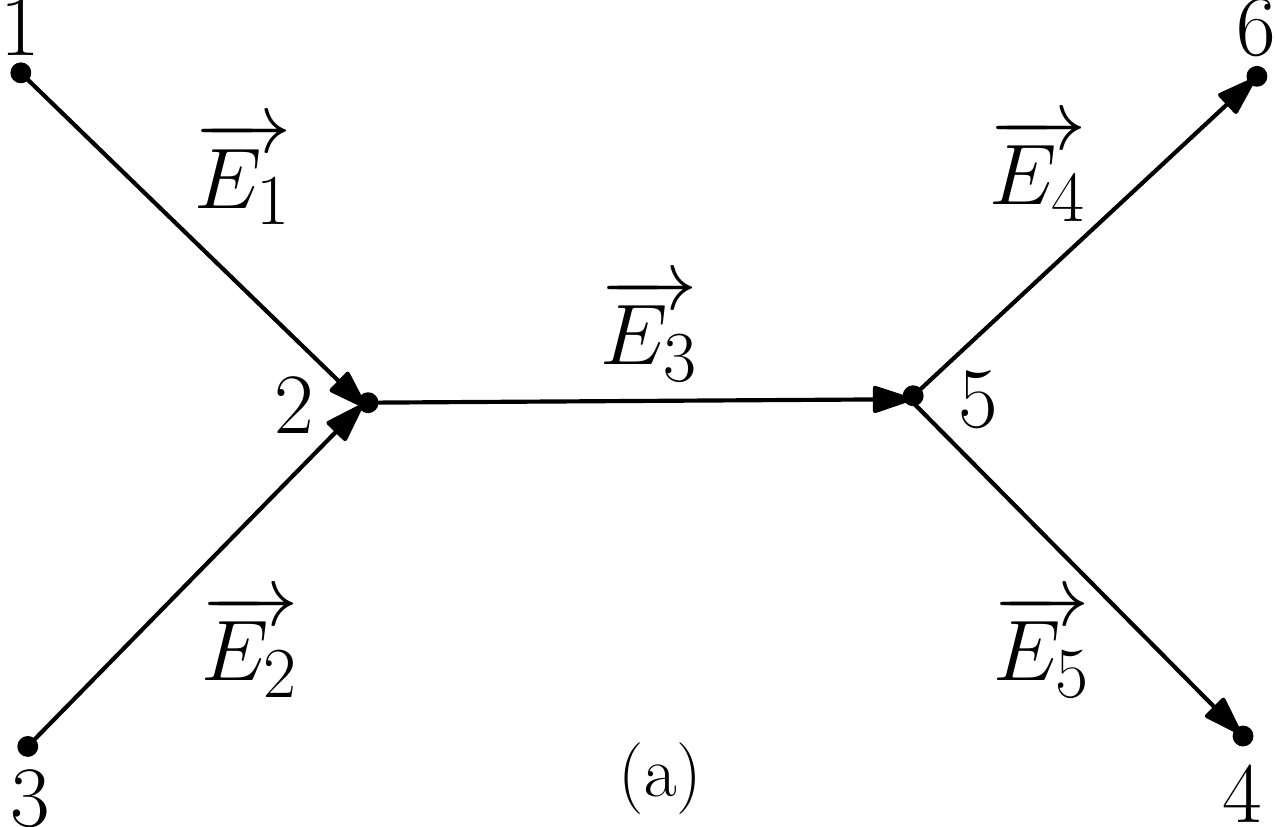}
\end{minipage}
\hspace{1in}
\begin{minipage}[t]{0.23\linewidth}
\centering
\includegraphics[width=1.6in, height=1.5in]{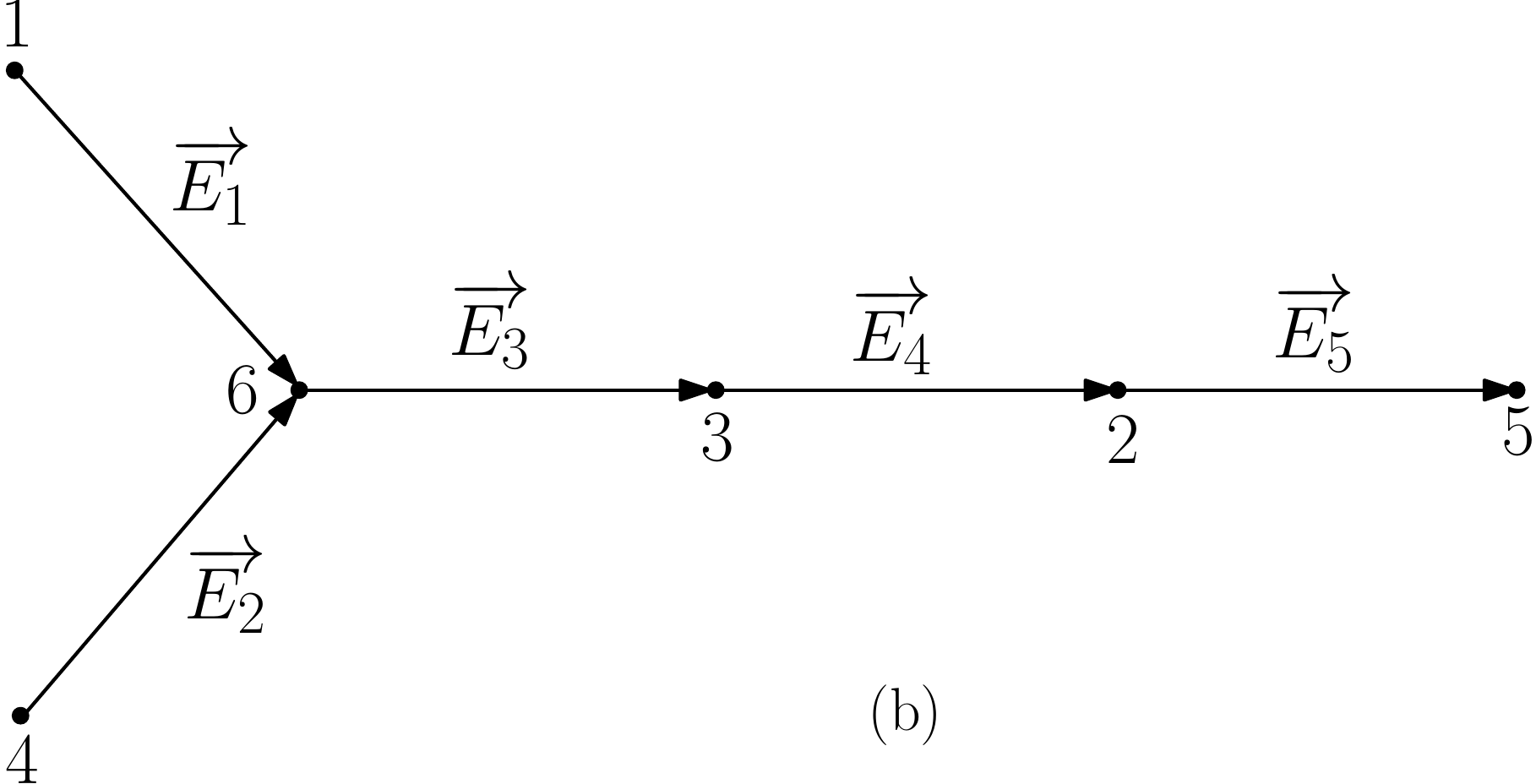}
\end{minipage}%
\hspace{1in}
\begin{minipage}[t]{0.23\linewidth}
\centering
\includegraphics[width=1.6in, height=1.5in]{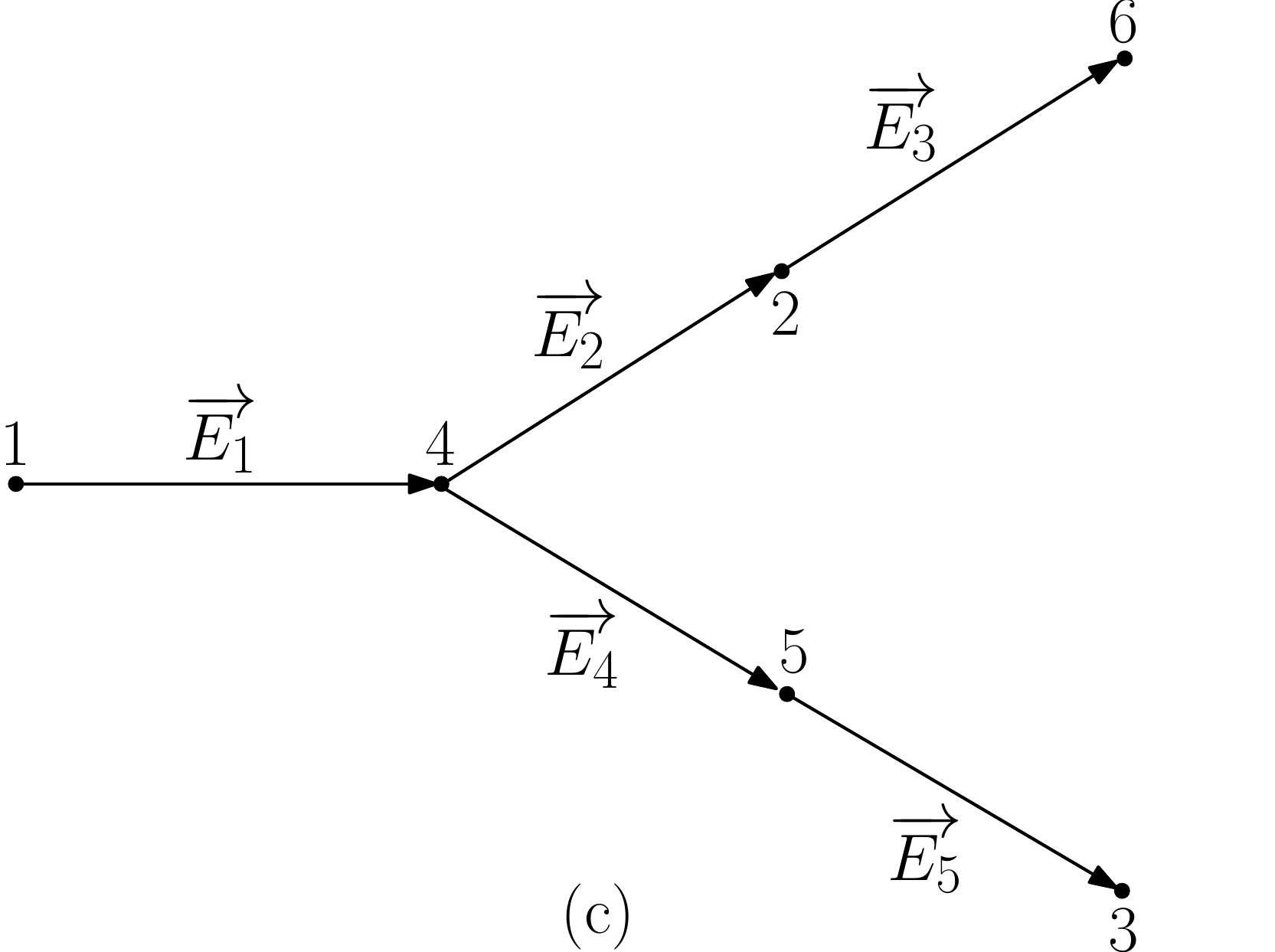}
\end{minipage}
\vspace{0.1in}

\[ \left[ {\begin{array}{*{20}c} 0 &  -1 &  0 &  0 &  0 \\ 0 &  -1 &  -1 & 0 & -1 \\ 0 &  1 & 1 &  1 & 0 \\ -1 &  0 &  -1 &  -1 &  0 \\ 0 &  0 &  0 &  -1 &  0 \\   \end{array}} \right],  \qquad\quad\, \left[ {\begin{array}{*{20}c} 
-1 & 0 &  -1 & -1 &  0 \\ -1 & 0 & -1 & -1 &  -1 \\ 1 &  -1 & 0 &  0 & 0 \\ 0 &  1 &  1 &  0 &  0 \\ 0 &  0 & -1 &  0 & 0 \\ \end{array}} \right], \qquad\quad\, \left[ {\begin{array}{*{20}c} 0 &  -1 &  0 &  1 &  0 \\ 0 &  0 &  0 & 0 &  1 \\ -1 &  0 & 0 &  -1 & -1 \\ 0 & 1 &  1 &  -1 &  0 \\ 0 &  -1 &  -1 &  0 &  0 \\  \end{array}} \right]\]
\vspace{0.8in}

\begin{minipage}[t]{0.23\linewidth}
\centering
\includegraphics[width=1.6in, height=1.5in]{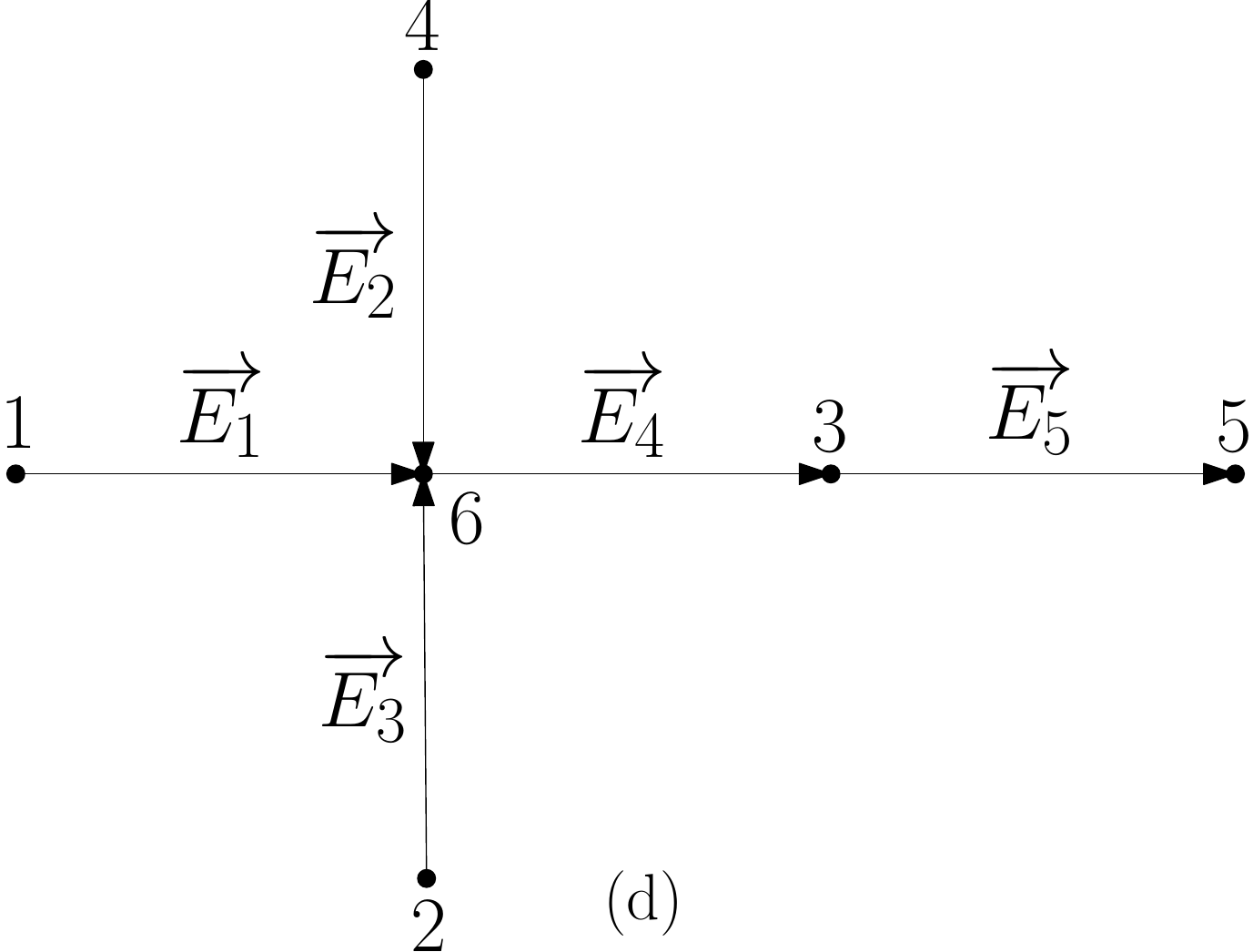}
\end{minipage}
\hspace{1in}
\begin{minipage}[t]{0.23\linewidth}
\centering
\includegraphics[width=1.6in, height=1.5in]{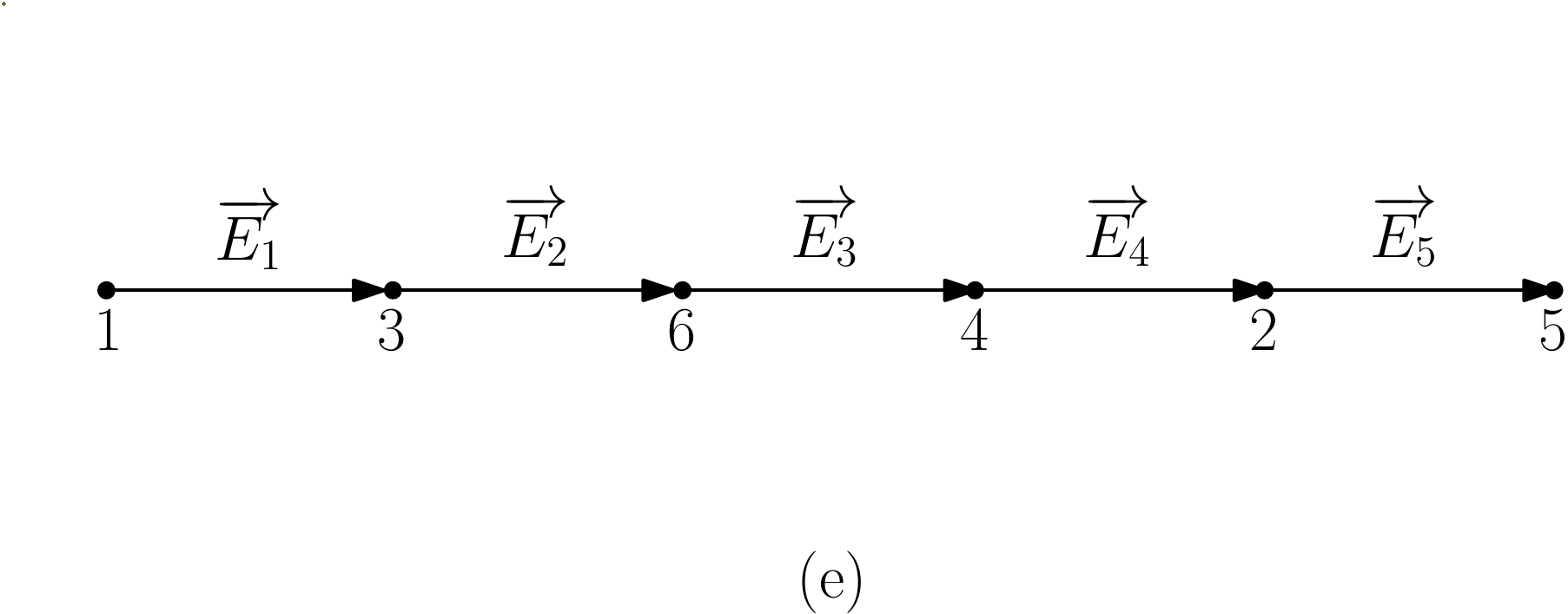}
\end{minipage}%
\hspace{1in}
\begin{minipage}[t]{0.23\linewidth}
\centering
\includegraphics[width=1.6in, height=1.5in]{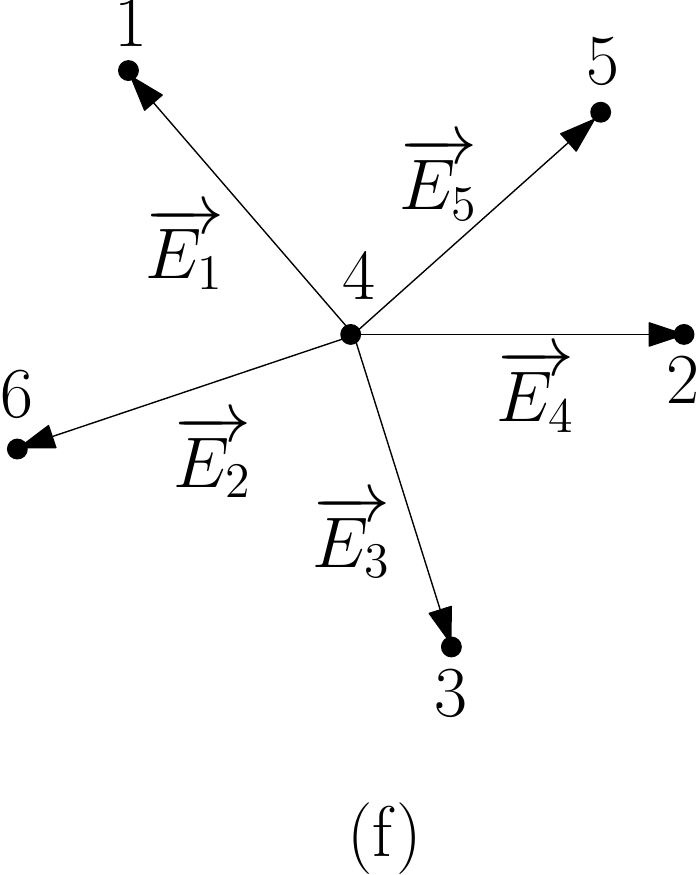}
\end{minipage}
\vspace{0.1in}

\[ \left[ {\begin{array}{*{20}c}  -1 &  0 &  1 & 0 & 0 \\ -1 &  0 & 0 & -1 &  -1 \\ -1 &  0 & 0 & -1 & 0 \\ 1 &  -1 &  0 & 0 & 0 \\ 0 &  1 &  0 &  0 &  0 \\  \end{array}} \right], \qquad\qquad\, \left[ {\begin{array}{*{20}c} 0 &  0 & 0 & -1 &  0 \\ -1 &  -1 & -1 & 0 &  0 \\ 1 &  1 & 1 &  1 & 1 \\ 0 &  -1 & -1 & -1 &  -1 \\ 0 & 1 &  0 &  0 &  0 \\ \end{array}} \right], \qquad\qquad \left[ {\begin{array}{*{20}c} 0 &  0 &  0 &  1 &  -1 \\ 1 &  0 &  0 & 0 &  -1 \\ 0 &  0 & 0 & 0 & -1 \\ 0 &  0 &  1 & 0 &  -1 \\ 0 &  1 &  0 &  0 &  -1 \\  \end{array}} \right]\]

\bigskip
\centering
\caption{The above 6 matrices are part of oriented transition matrices of continuous vertex maps on the oriented trees with 6 vertices right above them.  They are all similar to one another over the field $\mathcal F$.  The characteristic polynomials of their corresponding {\it unoriented transition matrices} are (a) $x^5-3x^4+x^3+x^2-3x+1$; (b) $x^5-x^4-3x^3-3x^2+x-1$; (c) $x^5-x^4-3x^3+x^2-x+1$; (d) $x^5-x^4-3x^3+x^2+x+1$; (e) $x^5-3x^4+x^3-3x^2-x-1$; and (f) $x^5-x^4-x^3-x^2-x-1$ respectively.}
\end{figure}

\begin{figure}[htbp]  
\centering
\includegraphics[width=3.0in, height=1.5in]{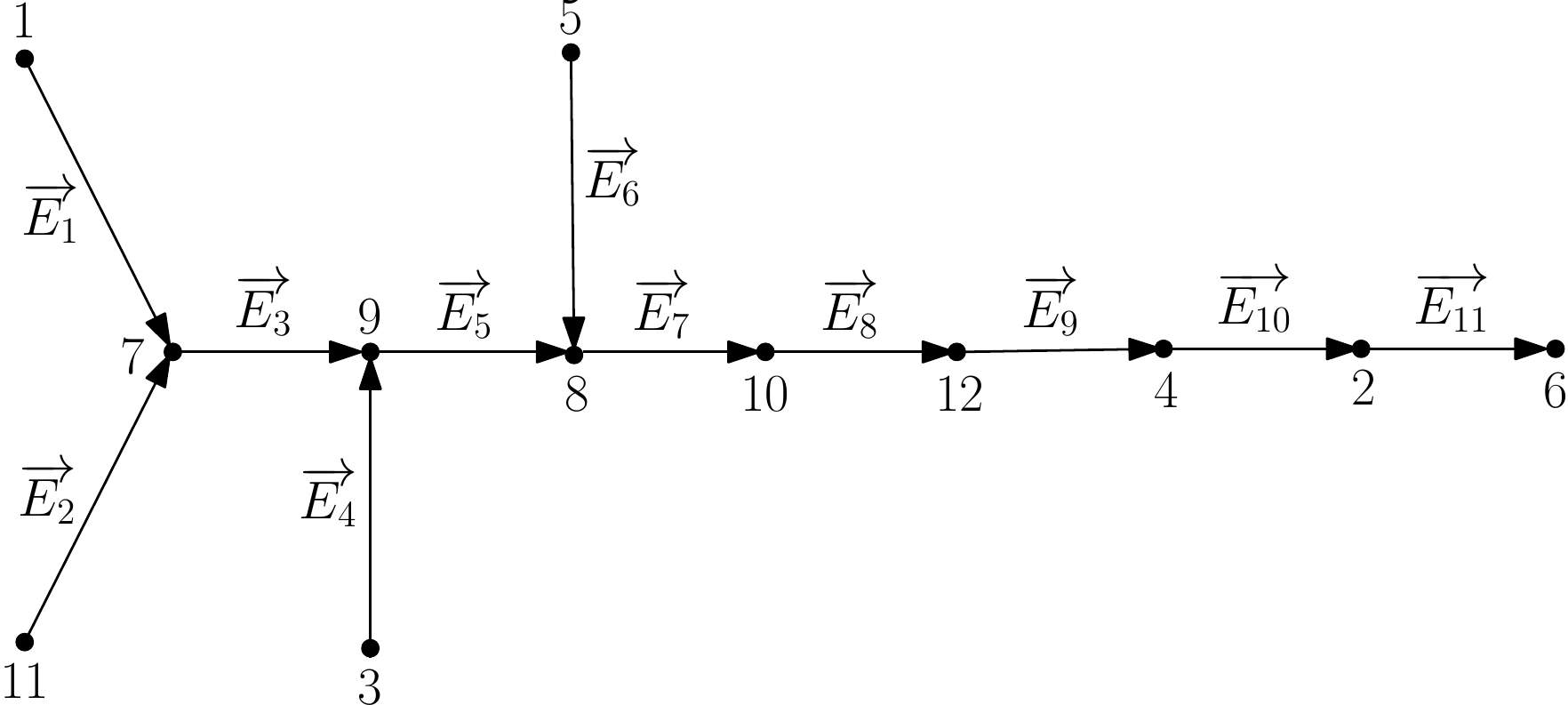}
\hspace{0.6in}
\includegraphics[width=3.0in, height=1.5in]{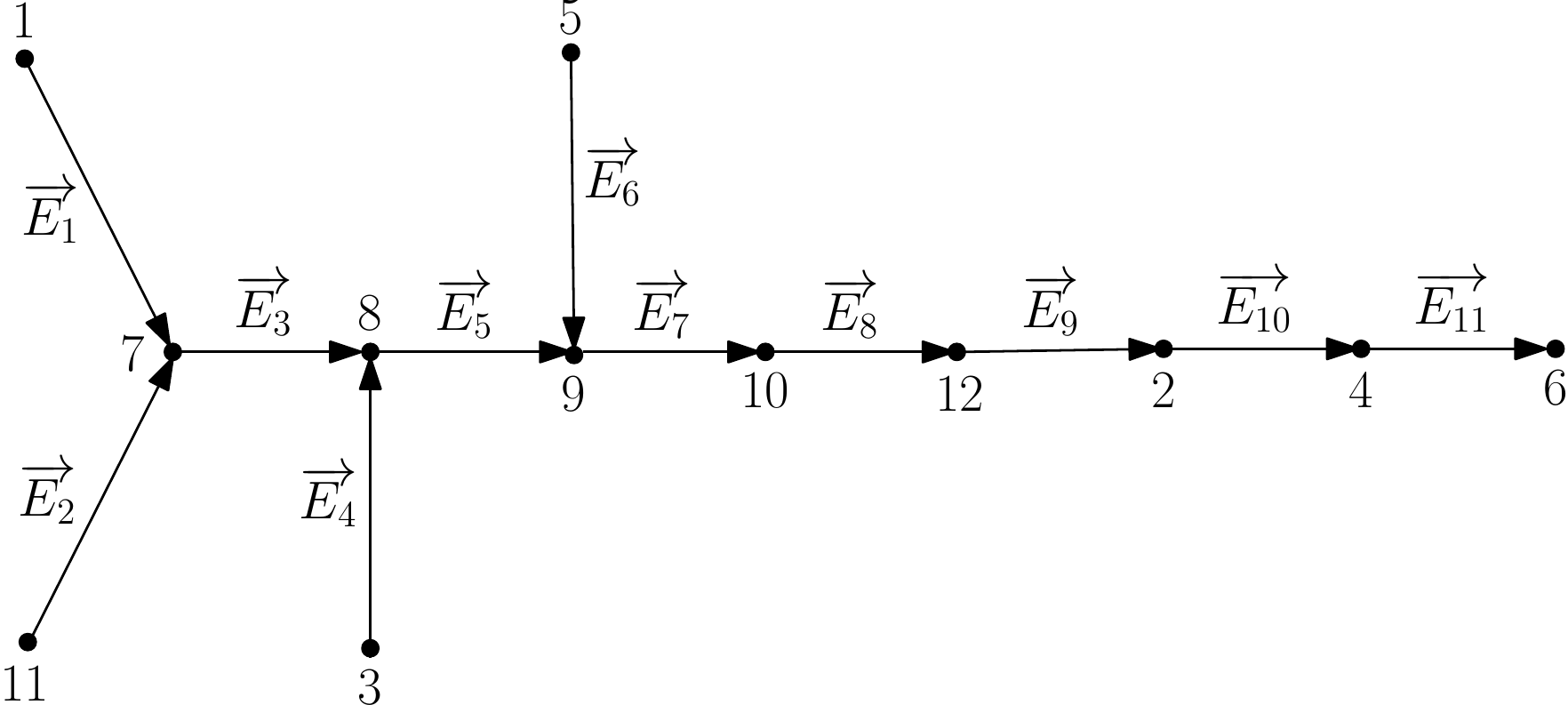}
\vspace{0.1in}
\scriptsize{
\[ \left[ {\begin{array}{*{20}c} 0 &  0 &  0 &  0 & 0 &  0 &  -1 &  -1 & -1 &  -1 &  0  \\ 0 &  0 &  0 &  0 & 0 &  0 &  -1 &  -1 & 0 &  0 &  0  \\ 0 &  0 &  0 &  0 & 0 &  0 &  1 &  0 & 0 &  0 &  0 \\ 0 &  0 &  0 &  0 & 0 &  0 &  0 &  -1 & -1 &  0 &  0 \\  0 &  0 &  0 &  0 & -1 &  0 &  -1 &  0 & 0 &  0 &  0 \\ 0 &  0 &  0 &  0 & -1 &  0 &  -1 &  -1 & -1 &  -1 &  -1 \\ 0 &  -1 &  -1 &  0 & 0 &  0 &  0 &  0 & 0 &  0 &  0 \\ -1 &  1 &  0 &  0 & 0 &  0 &  0 &  0 & 0 &  0 &  0 \\ 1 &  0 &  1 &  0 & 1 &  -1 &  0 &  0 & 0 &  0 &  0 \\ 0 &  0 &  0 &  -1 & -1 &  1 &  0 &  0 & 0 &  0 &  0 \\ 0 &  0 &  -1 &  1 & 0 &  0 &  0 &  0 & 0 &  0 &  0 \\    \end{array}} \right], \qquad  \left[ {\begin{array}{*{20}c} 0 &  0 &  0 &  0 & -1 &  0 &  -1 &  -1 & -1 &  0 &  0  \\ 0 &  0 &  0 &  0 & -1 &  0 &  -1 &  -1 & 0 &  0 &  0  \\ 0 &  0 &  0 &  0 & 1 &  0 &  0 &  0 & 0 &  0 &  0 \\ 0 &  0 &  0 &  0 & 0 &  0 &  -1 &  -1 & -1 &  -1 &  0 \\  0 &  0 &  0 &  0 & 0 &  0 &  1 &  0 & 0 &  0 &  0 \\ 0 &  0 &  0 &  0 & 0 &  0 &  0 &  -1 & -1 &  -1 &  -1 \\ 0 &  -1 &  -1 &  0 & -1 &  0 &  -1 &  0 & 0 &  0 &  0 \\ -1 &  1 &  0 &  0 & 0 &  0 &  0 &  0 & 0 &  0 &  0 \\ 1 &  0 &  1 &  -1 & 0 &  0 &  0 &  0 & 0 &  0 &  0 \\ 0 &  0 &  0 &  1 & 1 &  -1 &  0 &  0 & 0 &  0 &  0 \\ 0 &  0 &  -1 &  0 & -1 &  1 &  0 &  0 & 0 &  0 &  0 \\ \end{array}} \right]\] 
}
\centering
\caption{The characteristic polynomials of the above corresponding {\it unoriented transition matrices} are $x^{11}-x^{10}-7x^9+7x^8+7x^7-7x^6+3x^5-7x^4-x^3-x^2-3x+3$ and $x^{11}-x^{10}-9x^9+5x^8+25x^7-x^6-25x^5-11x^4+9x^3+9x^2-x-3$ respectively.}
\end{figure}%

\begin{figure}[htbp]  
\centering
\includegraphics[width=3.0in, height=1.5in]{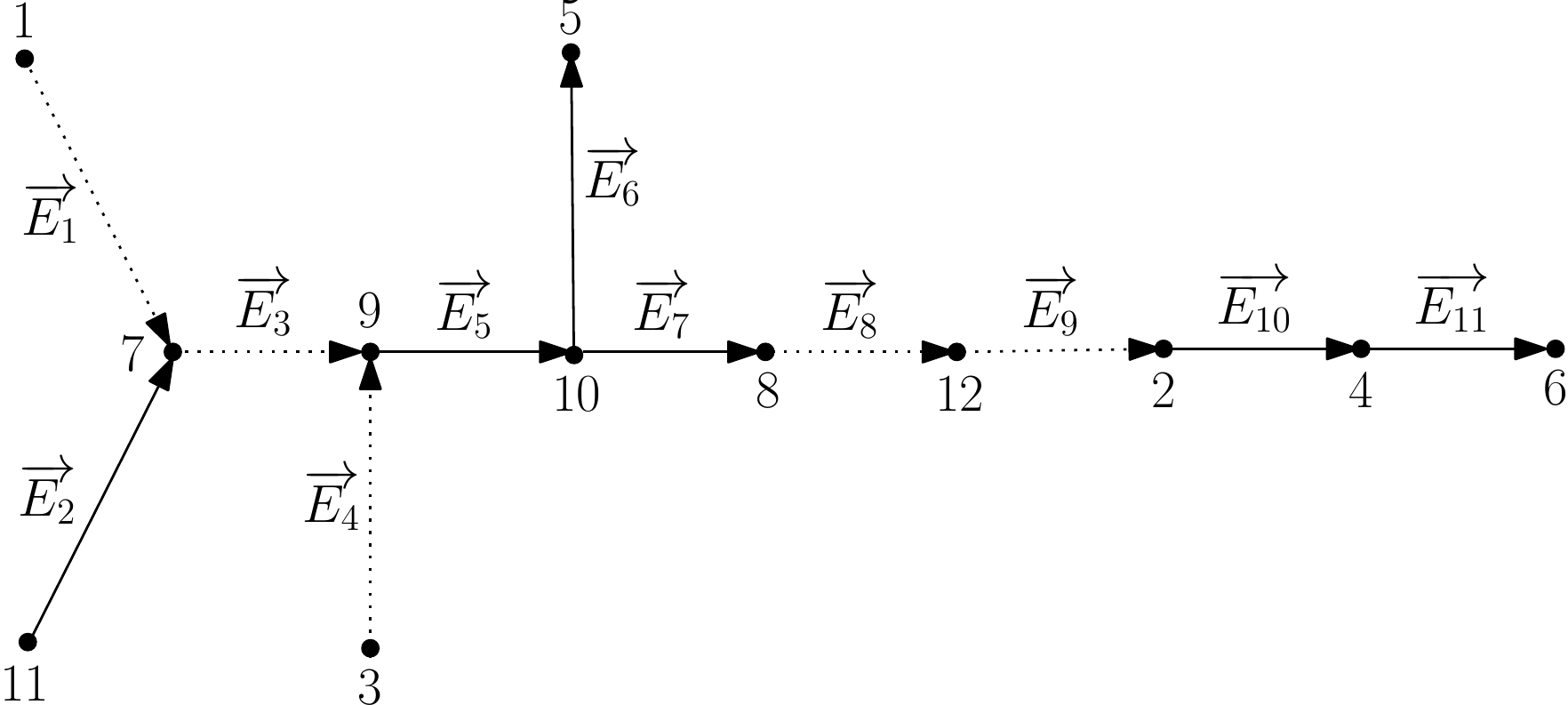}
\hspace{0.6in}
\includegraphics[width=3.0in, height=1.5in]{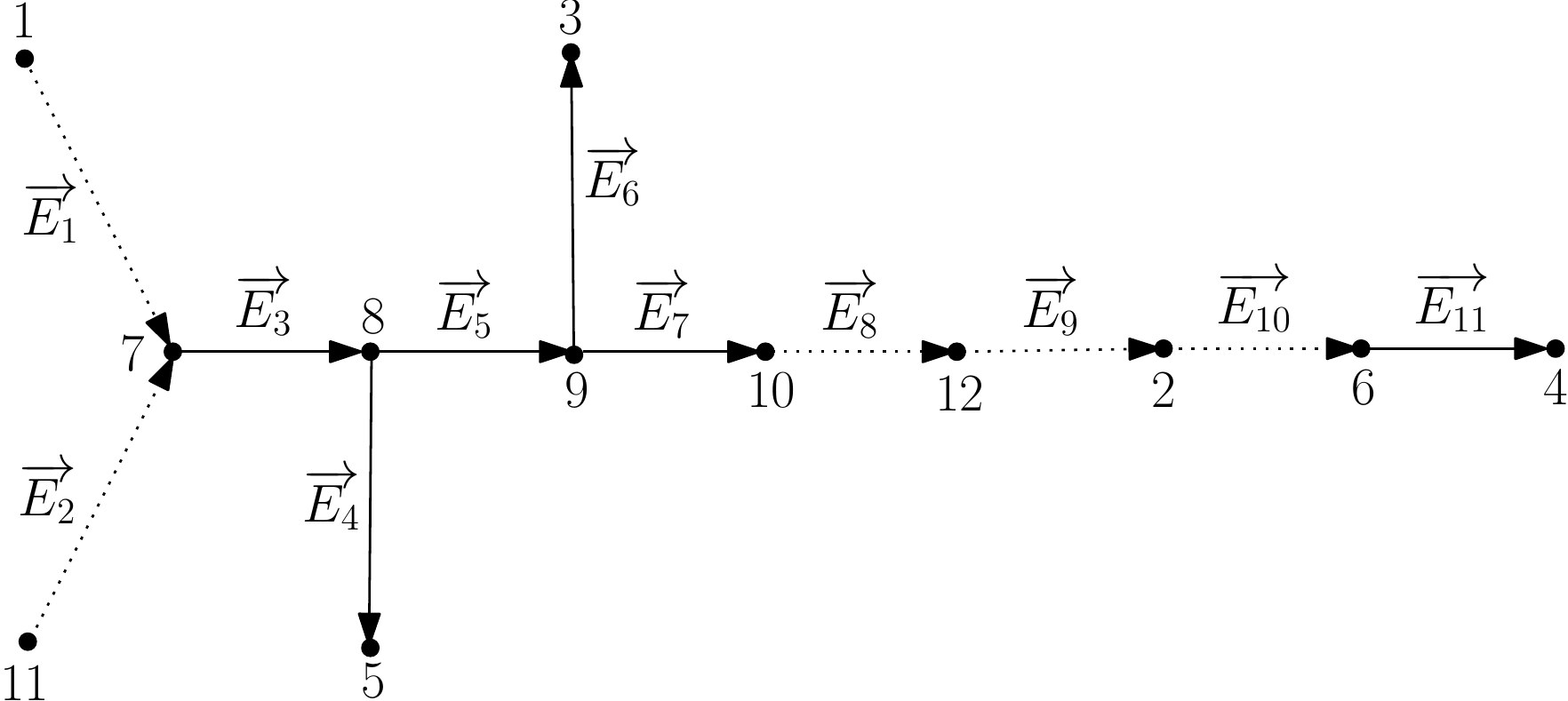}
\vspace{0.1in}
\footnotesize{   
\[ \left[ {\begin{array}{*{20}c} 0 &  0 &  0 &  0 & 0 &  0 &  0 &  -1 & -1 &  0 &  0  \\ 0 &  0 &  0 &  0 & 0 &  0 &  0 &  -1 & 0 &  0 &  0  \\ 0 &  0 &  0 &  0 & 0 &  0 &  -1 &  0 & 0 &  0 &  0 \\ 0 &  0 &  0 &  0 & 0 &  0 &  -1 &  -1 & -1 &  -1 &  0 \\  0 &  -1 &  -1 &  0 & -1 &  0 &  0 &  0 & 0 &  0 &  0 \\ 0 &  1 &  1 &  0 & 1 &  0 &  1 &  1 & 1 &  1 &  1 \\ 0 &  1 &  1 &  0 & 0 &  0 &  0 &  0 & 0 &  0 &  0 \\ -1 &  0 &  -1 &  0 & 0 &  0 &  0 &  0 & 0 &  0 &  0 \\ 1 &  0 &  1 &  -1 & 0 &  0 &  0 &  0 & 0 &  0 &  0 \\ 0 &  0 &  0 &  1 & 1 &  1 &  0 &  0 & 0 &  0 &  0 \\ 0 &  0 &  -1 &  0 & -1 &  -1 &  0 &  0 & 0 &  0 &  0 \\    \end{array}} \right], \,\,\,  \left[ {\begin{array}{*{20}c} 0 &  0 &  0 &  0 & -1 &  0 &  -1 &  -1 & -1 &  0 &  0  \\ 0 &  0 &  0 &  0 & -1 &  0 &  -1 &  -1 & 0 &  0 &  0  \\ 0 &  0 &  0 &  0 & 1 &  0 &  0 &  0 & 0 &  0 &  0 \\ 0 &  0 &  0 &  0 & 0 &  0 &  1 &  1 & 1 &  1 &  0 \\  0 &  0 &  0 &  0 & 0 &  0 &  1 &  0 & 0 &  0 &  0 \\ 0 &  0 &  0 &  0 & 0 &  0 &  0 &  1 & 1 &  1 &  1 \\ 0 &  -1 &  -1 &  0 & -1 &  0 &  -1 &  0 & 0 &  0 &  0 \\ -1 &  1 &  0 &  0 & 0 &  0 &  0 &  0 & 0 &  0 &  0 \\ 1 &  0 &  1 &  0 & 1 &  1 &  0 &  0 & 0 &  0 &  0 \\ 0 &  0 &  -1 &  0 & -1 &  -1 &  0 &  0 & 0 &  0 &  0 \\ 0 &  0 &  1 &  1 & 0 & 0 &  0 &  0 & 0 &  0 &  0 \\ \end{array}} \right]\] 
}
\centering
\caption{The characteristic polynomials of the above corresponding {\it unoriented transition matrices} are $x^{11}-x^{10}-7x^9+7x^8+13x^7-13x^6-7x^5+5x^4+x^3-x^2-x+1$ and $x^{11}-x^{10}-7x^9+3x^8+11x^7+5x^6+x^5-5x^4-5x^3-3x^2+x+1$ respectively.}
\end{figure}

\begin{figure}[htbp]  
\centering
\includegraphics[width=2.6in, height=1.5in]{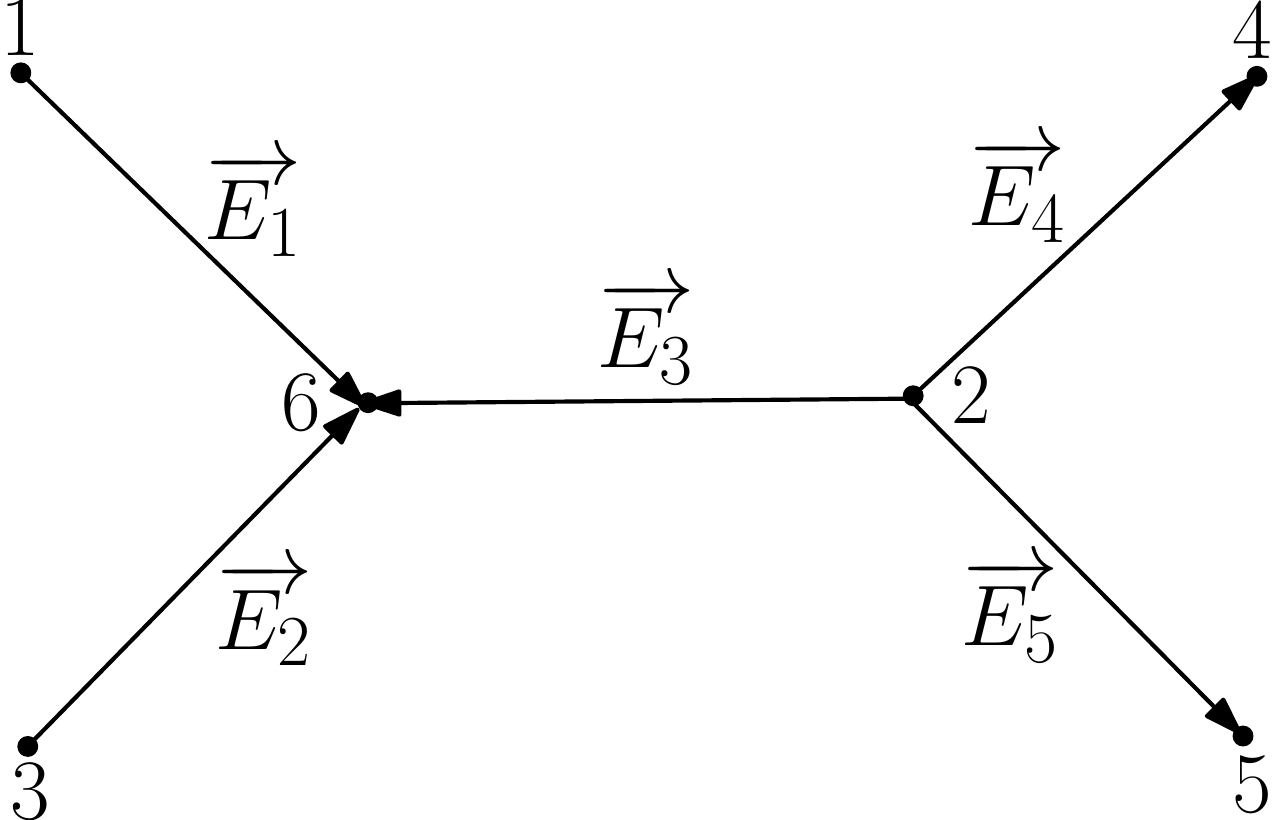}
\begin{equation*}
\left[ {\begin{array}{*{20}c} 1 &  0 &  -2 &  0 & 2  \\ 2 &  -1 &  -2 &  0 & 2  \\ 0 & 0 &  -1 &  0 & 2 \\ -2 &  0 &  2 &  1 & -2 \\  0 &  0 &  0 &  0 & 1  \\  \end{array}} \right] \left[ {\begin{array}{*{20}c} 1 &  0 &  1 &  0 & 0  \\ 1 &  0 &  1 &  1 & 0  \\ 1 &  1 &  0 &  0 & 0 \\ 0 &  1 &  1 &  0 & 1  \\  0 &  1 &  0 &  0 & 0  \\ \end{array}} \right] = \left[ {\begin{array}{*{20}c} -1 &  0 &  1 &  0 & 0  \\ -1 &  0 &  1 &  -1 & 0  \\ -1 &  1 &  0 &  0 & 0 \\ 0 &  1 &  -1 &  0 & 1  \\  0 &  1 &  0 &  0 & 0  \\  \end{array}} \right]   \end{equation*}
\centering
\caption{The characteristic polynomial of the above corresponding {\it unoriented transition matrices} is $x^5-x^4-3x^3-x^2+3x+1$.}
\end{figure}%

\pagebreak

\noindent
{\bf Acknowledgement}\\
This work was partially supported by the Ministry of Science and Technology of Taiwan.

\end{document}